\documentclass[a4paper,10pt]{amsart}
\usepackage{amsmath}
\usepackage{amsthm}
\usepackage{graphicx}
\usepackage{float}
\title{Properties of Bourbaki's Function}
\author{James McCollum}
\address{Division of Mathematics and Computer Science, University of Maine at Farmington, Farmington, ME 04938}
\email{james.mccollum@maine.edu}
\keywords{continuous, nowhere differentiable, fractal, box dimension, Hausdorff dimension}
\subjclass{26A30,37C45}
\date{September 2010}
\begin{document}

\maketitle

\begin{abstract}
We examine Bourbaki's function, an easily-constructed continuous but nowhere-differentiable function, and explore properties including functional identities, the antiderivative, and the box and Hausdorff dimensions of the graph.
\end{abstract}

\section{Introduction}
While Bernard Bolzano \cite{Jarnik81} introduced one of the earliest examples of a continuous, nowhere-differentiable function, his example is only one of countless similar functions, many of which are defined in less complex ways.  One such function is found in Nicolas Bourbaki's \emph{Elements of Mathematics---Functions of a Real Variable} \cite{Bourbaki04}: This function, which we call Bourbaki's Function, is defined by a few inductive rules, and its simple self-similar structure allows for abundant and relatively easy analysis.

Okamoto \cite{Okamoto05} defines Bourbaki's Function $f_i$ for any iteration $i\geq0$ over $[0,1]$ as follows: $f_0(x)=x$ for all $x \in [0,1]$, every $f_i$ is continuous on $[0,1]$, every $f_i$ is affine in each subinterval $[k/3^i,(k+1)/3^i]$ where $k\in\{0,1,2,\ldots,3^i-1\}$, and
\begin{flalign}
&f_{i+1}\left(\frac{k}{3^i}\right)=f_i\left(\frac{k}{3^i}\right),\\
&f_{i+1}\left(\frac{3k+1}{3^{i+1}}\right)=f_i\left(\frac{k}{3^i}\right)+\frac{2}{3}\left[f_i\left(\frac{k+1}{3^i}\right)-f_i\left(\frac{k}{3^i}\right)\right],\\
&f_{i+1}\left(\frac{3k+2}{3^{i+1}}\right)=f_i\left(\frac{k}{3^i}\right)+\frac{1}{3}\left[f_i\left(\frac{k+1}{3^i}\right)-f_i\left(\frac{k}{3^i}\right)\right],\\
&f_{i+1}\left(\frac{k+1}{3^i}\right)=f_i\left(\frac{k+1}{3^i}\right).
\end{flalign}
Figs. 1 and 2 illustrate the construction of the graph of $f$.  Okamoto \cite{Okamoto05} has shown, using the above equations, that the function $$f(x)=\displaystyle\lim_{i\to\infty}f_i(x)$$ is continuous and nowhere differentiable.  We can observe from these equations that $f_i(x)=f(x)$ for any $x\in[0,1]$ that can be expressed as some integral multiple of $1/3^i$.  The values between each of these inputs, however, may change with new iterations.

\begin{figure}[t]
\begin{center}
\includegraphics[width=.75\textwidth]{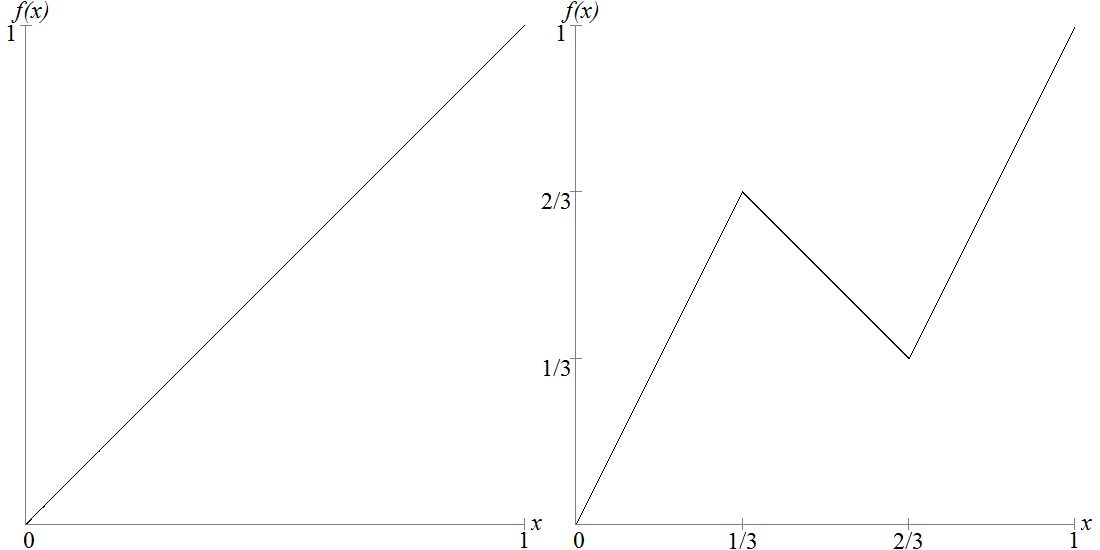}
\caption{Graphs of $f_0$ and $f_1$. Note the steps taken in constructing $f_1$ from $f_0$.}
\end{center}
\end{figure}

\begin{figure}[h!]
\begin{center}
\includegraphics[width=.75\textwidth]{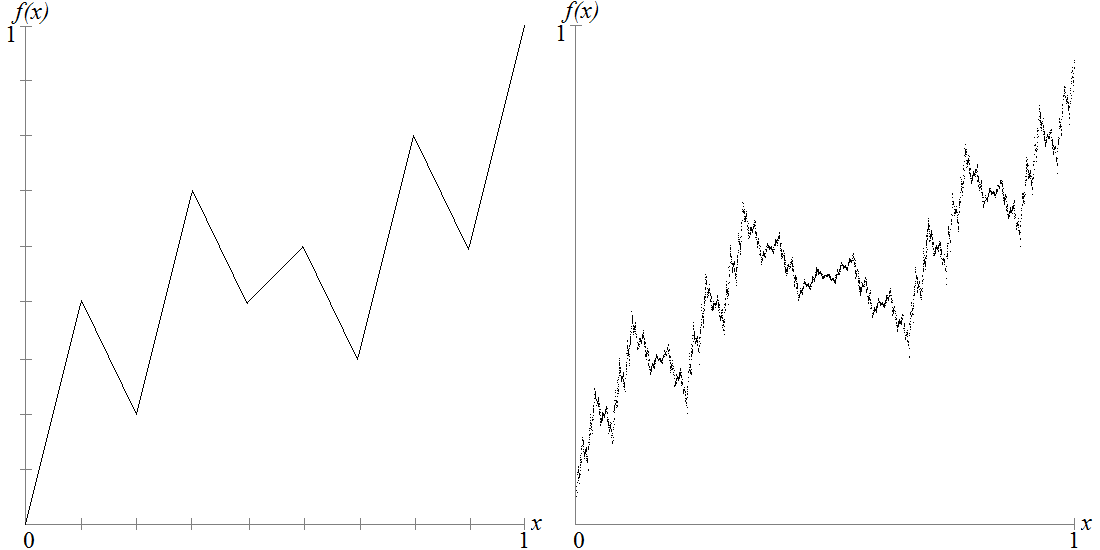}
\caption{Graph of $f_2$ and an approximate graph of $f$.  Again, note the steps taken in constructing $f_2$ from $f_1$.}
\end{center}
\end{figure}

Formulas (1)--(4) can evaluate $f(x)$ easily for integral multiples of $1/3^i$: Consider the ternary expansion $x=0.x_1x_2\ldots x_i$, where $x_1,x_2,\ldots,x_i \in \{0,1,2\}$. (If $x=1$, we say equivalently that $x=0.222\ldots$) The value of $i$ coincides with the first iteration $i$ for which $f_i(x)=f(x)$.  We can rewrite the ternary expansion as
\begin{equation*}
x=\sum_{j=1}^i\frac{x_j}{3^j}.
\end{equation*}
Since the sum always will have a denominator of $1/3^i$, we can evaluate $f(x)$ using an $f_i$ that covers intervals of length at least $1/3^i$---which, by definition, is $f_h$ for $h\geq i$.  To evaluate $f(x)$, we start with $k=0$ and $i=1$, and we apply formula (1), (2), or (3) depending on the value of each $x_j$: If $x_j=0$, we use (1); if $x_j=1$, we use (2); and if $x_j=2$, we use (3).  In every case, we must keep track of the values $k/3^i$, $(3k+1)/3^{i+1}$, $(3k+2)/3^{i+1}$, $(3k+1)/3^{i}$, and $(3k+2)/3^{i}$.

To evaluate $f(x)$ for other $x$ values---when $x$ is an irrational number or any rational number in $[0,1]$ whose denominator is not a power of three---we must consider a non-terminating ternary expansion:
\begin{equation*}
x=\sum_{j=1}^\infty\frac{x_j}{3^j}.
\end{equation*}
In other words, to evaluate $f(x)$ for such an $x$, we would have to apply formulas (1), (2), and (3) indefinitely.  While we certainly could obtain a fair approximation with enough iterations, keeping track of certain values---in particular, $f_i((k+1)/3^i)-f_i(k/3^i)$---would grow more difficult with each iteration.

In this paper, we will use the self-similarity of $f$ to find ``shortcuts'' to evaluating $f(x)$ for such values of $x$, and we will examine how this self-similarity concerns the fractal nature of the graph of $f$, the antiderivative $F$ of $f$, and the properties of the graph of $F$.  More specifically, in Section 2, we will prove that the graph of $f$ possesses rotational symmetry about the point $(1/2,1/2)$---in other words, 
\begin{equation*}
f(1-x)=1-f(x) \textrm{\ for all\ }x\in[0,1].
\end{equation*}
In the same section, we will use the function's self-similar properties to infer three basic identities that evaluate $f(x/3^i)$, $f([2-x]/3^i)$, and $f([2+x]/3^i)$ in terms of $f(x)$.  In Section 3, we will use these general identites to evaluate $f(x)$ for specific sets of numbers that have some form other than $x=k/3^i$.  In Section 4, we will use the rotational symmetry of the graph of $f$ to prove that 
\begin{equation*}
\int_x^{1-x}f(t)\,dt=1/2-x.
\end{equation*}
In the same section, we will derive three other identities for the area under the graph of $f$, and we will use these identities iteratively to construct a graph of $F$.  In Section 5, we will do for $F$ what we did for $f$ in Section 3.  Finally, in Section 6, we will show that the graph of $f$ has box-counting and Hausdorff dimensions equal to $\log_3 5$, and we will show that the arc length of the graph of $F$ is bounded below by $\displaystyle\frac{\sqrt{5}}{2}$ and above by $\displaystyle\frac{3}{2}$.

\section{Functional Identities}
We observe from the graphs that each $f_i$ possesses rotational symmetry about its center, which implies a useful identity for $f$:
\newtheorem{theorem}{Theorem}
\begin{theorem}
For all $x \in [0,1]$, $f(1-x)=1-f(x)$.
\begin{proof}
By definition, $f_0(x)=x$ for all $x \in [0,1]$, so clearly this is true for $f_0$.  To prove this for $f$ itself, however, we will take advantage of the function's inductive nature and consider only the points where $x$ can be expressed as a rational number in $[0,1]$ whose denominator takes the form $3^i$.

We consider $i=1$ for a base case. We must show, then, that $f_1(1-x)=1-f(x)$ for $x \in \{0,1/3,2/3,1\}$.  Evaluating the function at these points, we get
\begin{flalign*}
f_1(1-0)&=f_1(1)=1=1-0=1-f_1(0)\\
f_1\left(1-\frac{1}{3}\right)&=f_1\left(\frac{2}{3}\right)=\frac{1}{3}=1-\frac{2}{3}=1-f_1\left(\frac{1}{3}\right)\\
f_1\left(1-\frac{2}{3}\right)&=f_1\left(\frac{1}{3}\right)=\frac{2}{3}=1-\frac{1}{3}=1-f_1\left(\frac{2}{3}\right)\\
f_1(1-1)&=f_1(0)=0=1-1=1-f_1(1)
\end{flalign*}
So for $i=1$, $f_i(1-x)=1-f_i(x)$ for all $x \in \{0,1/3^i,2/3^i,\ldots,1\}$.

From here we can make our inductive hypothesis: For some $j\geq1$, $f_j(1-x)=1-f_j(x)$, where $x \in \{0,1/3^j,2/3^j,\ldots,1\}$.  Now, we can deduce from (1) and (4) that if $f_j(1-x)=1-f_j(x)$ and $f_{j+1}(x)=f_j(x)$ for $x=k/3^j$, then $f_{j+1}(1-x)=1-f_{j+1}(x)$, since both $x,(1-x)\in\{0,1/3^j,2/3^j,\ldots,1\}$.  To apply (2) and (3), we first let $k'=3^j-1-k$. This means that $k' \in \{0,1,\ldots,3^j-1\}$ and that $1-k/3^j=(k'+1)/3^j$.  Then
\begin{flalign*}
f_{j+1}\left(1-\frac{3k+1}{3^{j+1}}\right) &=f_{j+1}\left(\frac{3^{j+1}-3k-1}{3^{j+1}}\right)\\
&=f_{j+1}\left(\frac{(3^{j+1}-3k-3)+2}{3^{j+1}}\right)\\
&=f_{j+1}\left(\frac{3k'+2}{3^j}\right)\\
&=f_j \left(\frac{k'}{3^j} \right)+\frac{1}{3}\left[f_j\left(\frac{k'+1}{3^j}\right)-f_j\left(\frac{k'}{3^j}\right)\right].
\end{flalign*}
Using our inductive hypothesis, we make appropriate substitutions for these iteration-$j$ functions to get
\begin{flalign*}
f_{j+1}\left(1-\frac{3k+1}{3^{j+1}}\right)&=1-f_j\left(1-\frac{k'}{3^j} \right)+\frac{1}{3}\left[1-f_j\left(1-\frac{k'+1}{3^j}\right)\right]\\
&\hspace{1 pc}-\frac{1}{3}\left[1-f_j\left(1-\frac{k'}{3^j}\right)\right]\\
&=1-f_j\left(\frac{k+1}{3^j}\right)+\frac{1}{3}f_j\left(\frac{k+1}{3^j}\right)-\frac{1}{3}f_j\left(\frac{k}{3^j}\right)\\
&=1-\frac{2}{3}f_j\left(\frac{k+1}{3^j}\right)+\frac{2}{3}f_j\left(\frac{k}{3^j}\right)-f_j\left(\frac{k}{3^j}\right)\\
&=1-\left(f_j\left(\frac{k}{3^j}\right)+\frac{2}{3}\left[f_j\left(\frac{k+1}{3^j}\right)-f_j\left(\frac{k}{3^j}\right)\right]\right)\\
&=1-f_{j+1}\left(\frac{3k+1}{3^{j+1}}\right).
\end{flalign*}
Working the same substitutions through for $f_{j+1}(1-(3k+2)/3^{j+1})$ will give us $1-f_{j+1}((3k+2)/3^{j+1})$. 

Therefore, by induction, $f_i(1-x)=1-f_i(x)$ for all $i\geq1$ and all $x\in\{0,1/3^i,2/3^i,\ldots,1\}$.  As $i$ approaches infinity, the interval $1/3^i$ between each $x\in\{0,1/3^i,2/3^i,\ldots,1\}$ approaches zero, and since this set is dense in $[0,1]$, the limit $f(x)$ satisfies $f(1-x)=1-f(x)$ for all $x \in [0,1]$.
\end{proof}
\end{theorem}
This identity proves helpful in evaluating $f(x)$ for values of $x$ whose ternary expansions do not terminate---used in conjunction with the three propositions of this section, this identity makes it possible to evaluate $f(x)$ for values of the form $1/(3^i+1)$, for instance.  We also can use it to show that $f(1/2)=1/2$:
\begin{flalign*}
f\left(\frac{1}{2}\right)&=1-f\left(\frac{1}{2}\right)\\
2f\left(\frac{1}{2}\right)&=1\\
f\left(\frac{1}{2}\right)&=\frac{1}{2}
\end{flalign*}
Katsuura \cite{Katsuura91} defines the contraction mappings $w_n: X \mapsto X$, where $n \in \{1,2,3\}$ and $X=[0,1] \times [0,1]$, as follows: For all $(x,y)\in X$,
\begin{flalign}
w_1(x,y)&=\left(\frac{x}{3},\frac{2y}{3}\right)\\
w_2(x,y)&=\left(\frac{2-x}{3},\frac{1+y}{3}\right)\\
w_3(x,y)&=\left(\frac{2+x}{3},\frac{1+2y}{3}\right)
\end{flalign}
Seperately applying mappings (5), (6), and (7) to the line $y=x$ (the graph of $f_0$) produces the graph of $f_1$; applying the same mappings to the graph of $f_1$ gives the graph of $f_2$; and so on.  More generally, if $\Gamma_i$ is the graph of $f_i$, then 
\begin{equation*}
\Gamma_{i+1}=w_1(\Gamma_i)\cup w_2(\Gamma_i)\cup w_3(\Gamma_i).
\end{equation*}
Since $f=\displaystyle\lim_{i\to\infty}f_i$, we can say that $\Gamma=\displaystyle\lim_{i\to\infty}\Gamma_{i+1}=\displaystyle\lim_{i\to\infty}w_1(\Gamma_i)\cup w_2(\Gamma_i)\cup w_3(\Gamma_i)$.  So $\Gamma=w_1(\Gamma)\cup w_2(\Gamma)\cup w_3(\Gamma)$ is the unique invariant set for the iterated function system (IFS) given by $w_1$, $w_2$, and $w_3$ (see \cite{Katsuura91}).  Since $w_1(\Gamma_i)=\Gamma_{i+1}$ on $[0,1/3]$, $w_2(\Gamma_i)=\Gamma_{i+1}$ on $[1/3,2/3]$, and $w_3(\Gamma_i)=\Gamma_{i+1}$ on $[2/3,1]$, we are able to prove three more identities:
\newtheorem{proposition}{Proposition}
\begin{proposition}
For all $x \in [0,1]$ and $i\geq 0$, $\displaystyle f\left(\left(\frac{1}{3}\right)^ix\right)=\left(\frac{2}{3}\right)^if(x)$.
\begin{proof}
Let $x\in[0,1]$.  Then $(x, f(x))\in\Gamma$.  If $i=0$, our result is obvious.  If $i>0$, then
\begin{equation*}
\underbrace{w_1\circ w_1\circ \dots \circ w_1}_{i-1}\circ w_1(x,f(x))=\left(\left(\frac{1}{3}\right)^ix,\left(\frac{2}{3}\right)^if(x)\right)
\end{equation*}
from the definition of $w_1$.  And since $w_1^i(\Gamma)\subseteq\Gamma$, where $w_1^i(\Gamma)$ denotes $i$ applications of $w_1$ on $\Gamma$, we have $\displaystyle f\left(\left(\frac{1}{3}\right)^ix\right)=\left(\frac{2}{3}\right)^if(x)$.
\end{proof}
\end{proposition}

\begin{proposition}
For all $x \in [0,1]$ and $i>0$, 
\begin{equation*}
f\left(\frac{2-x}{3^i}\right)=\frac{2^{i-1}}{3^i}[1+f(x)].
\end{equation*}
\begin{proof}
Let $x\in[0,1]$.  Then $(x, f(x))\in\Gamma$.  If $i>0$, then
\begin{flalign*}
\underbrace{w_1\circ w_1\circ \dots \circ w_1}_{i-1}\circ w_2(x,y)&=\left(\frac{([2-x]/3)}{3^{i-1}},\frac{2^{i-1}([1+y]/3)}{3^{i-1}}\right)\\
&=\left(\frac{2-x}{3^i},\frac{2^{i-1}(1+y)}{3^i}\right)
\end{flalign*}
from the definition of $w_1$.  And since $w_1(\Gamma)\subseteq\Gamma$ and $w_2^i(\Gamma)\subseteq\Gamma$, we have $\displaystyle f\left(\frac{2-x}{3^i}\right)=\frac{2^{i-1}}{3^i}[1+f(x)]$.
\end{proof}
\end{proposition}
\begin{proposition}
For all $x \in [0,1]$ and $i>0$, 
\begin{equation*}
f\left(\frac{2+x}{3^i}\right)=\left(\frac{2}{3}\right)^if(x)+\frac{2^{i-1}}{3^i}.
\end{equation*}
\begin{proof}
This can be proven in the same manner as Proposition 2 if we replace the $w_2$ with $w_3$.
\end{proof}
\end{proposition}

\section{Function Values}
As $f$ has no explicit formula, we must take advantage of its self-similar structure to evaluate $f(x)$ for nearly all values of $x$ (see Fig. 3).  The four identities we have proven already will help.

\begin{figure}[htb]
\begin{center}
\includegraphics[width=1.0\textwidth]{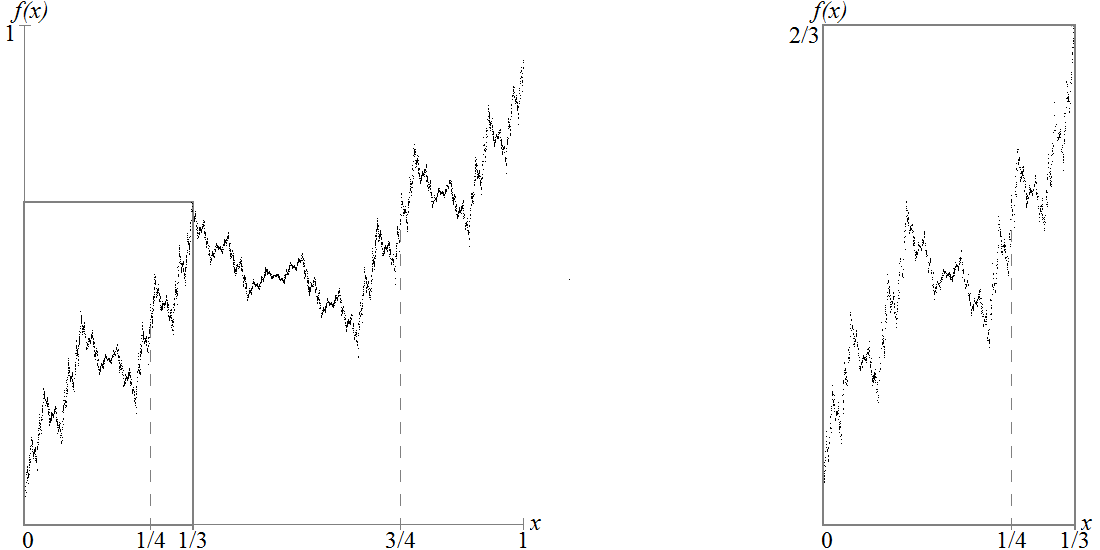}
\caption{As $f$ has no explicit formula, the self-similarity of its graph is key in determining its values at different points.}
\end{center}
\end{figure}

\begin{theorem}
For all $j>i>0$,
\renewcommand{\labelenumi}{(\roman{enumi})}
\begin{enumerate}
\item $\displaystyle f\left(\frac{1}{3^i+1}\right)=\frac{2^i}{3^i+2^i}$,\\
\item $\displaystyle f\left(\frac{1}{3^i-1}\right)=\frac{2^i}{3^i+2^{i-1}}$,\\
\item $\displaystyle f\left(\frac{2}{3^i+1}\right)=\frac{2^{i-1}}{3^i-2^{i-1}}$,\\
\item $\displaystyle f\left(\frac{2}{3^i-1}\right)=\frac{2^{i-1}}{3^i-2^i}$,\\
\item $\displaystyle f\left(\frac{1}{3^j+3^i}\right)=\left(\frac{2}{3}\right)^i\left(\frac{2^{j-i}}{3^{j-i}+2^{j-i}}\right)$, and\\
\item $\displaystyle f\left(\frac{1}{3^j-3^i}\right)=\left(\frac{2}{3}\right)^i\left(\frac{2^{j-i}}{3^{j-i}+2^{j-i-1}}\right)$.
\end{enumerate}
\begin{proof} Let $j>i>0$.\\

(i) Clearly
\begin{equation*}
1-\frac{1}{3^i+1}=\frac{3^i}{3^i+1}
\end{equation*}
and using Theorem 1 and Proposition 1, it follows that
\begin{flalign*}
f\left(\frac{1}{3^i+1}\right)&=f\left(\frac{1}{3^i}\left[\frac{3^i}{3^i+1}\right]\right)\\
&=f\left(\frac{1}{3^i}\left[1-\frac{1}{3^i+1}\right]\right)\\
&=\left(\frac{2}{3}\right)^if\left(1-\frac{1}{3^i+1}\right)\\
&=\left(\frac{2}{3}\right)^i\left[1-f\left(\frac{1}{3^i+1}\right)\right]\\
&=\left(\frac{2}{3}\right)^i-\left(\frac{2}{3}\right)^if\left(\frac{1}{3^i+1}\right).
\end{flalign*}
So we have
\begin{equation*}
\left[1+\left(\frac{2}{3}\right)^i\right]f\left(\frac{1}{3^i+1}\right)=\left(\frac{2}{3}\right)^i
\end{equation*}
and thus,
\begin{equation*}
f\left(\frac{1}{3^i+1}\right)=\frac{2^i}{3^i+2^i},\textrm{\ for\ }i>0.\\
\end{equation*}

(ii) Likewise, we know that for $i>0$,
\begin{equation*}
\frac{1}{3^i-1}=1-\frac{3^i-2}{3^i-1}
\end{equation*}
and
\begin{equation*}
\frac{2-(3^i-2)/(3^i-1)}{3^i}=\frac{1}{3^i-1}.
\end{equation*}
The next steps are almost identical to those from the previous proof, so we will omit them.  Making appropriate substitutions, applying the function to both sides, and using Theorem 1 and Proposition 2 gives us
\begin{equation*}
f\left(\frac{1}{3^i-1}\right)=\frac{2^i}{3^i+2^{i-1}},\textrm{\ for\ }i>0.\\
\end{equation*}

(iii) For $i>0$, we have

\begin{equation*}
\frac{2-(2/(3^i+1))}{3^i}=\frac{2}{3^i+1}.
\end{equation*}
Applying the function to both sides and using Proposition 2 gives us
\begin{equation*}
f\left(\frac{2}{3^i+1}\right)=\frac{2^{i-1}}{3^i-2^{i-1}},\textrm{\ for\ }i>0.\\
\end{equation*}

(iv) Similarly,
\begin{equation*}
\frac{2+(2/(3^i-1))}{3^i}=\frac{2}{3^i-1},
\end{equation*}
and by applying $f$ to both sides and using Proposition 3, we get
\begin{equation*}
f\left(\frac{2}{3^i-1}\right)=\frac{2^{i-1}}{3^i-2^i},\textrm{\ for\ }i>0.\\
\end{equation*}

(v) Now, for $j>i>0$, we know that
\begin{equation*}
\frac{1}{3^j+3^i}=\left(\frac{1}{3^i}\right)\left(\frac{1}{3^{j-i}+1}\right)
\end{equation*}
and since $j>i>0$ implies that $j-i>0$, we can apply Proposition 1 and the identity $f(1/(3^i+1))=2^i/(3^i+2^i)$ (which we proved in part (i) of this theorem) to obtain
 \begin{flalign*}
f\left(\frac{1}{3^j+3^i}\right)&=f\left(\left(\frac{1}{3^i}\right)\left(\frac{1}{3^{j-i}+1}\right)\right)\\
&=\left(\frac{2}{3}\right)^if\left(\frac{1}{3^{j-i}+1}\right)\\
&=\left(\frac{2}{3}\right)^i\left(\frac{2^{j-i}}{3^{j-i}+2^{j-i}}\right),\textrm{\ for\ }j>i>0.\\
\end{flalign*}

(vi) Given $j>i>0$, we also know that
\begin{equation*}
\frac{1}{3^j-3^i}=\left(\frac{1}{3^i}\right)\left(\frac{1}{3^{j-i}-1}\right)
\end{equation*}
and by applying Proposition 1 and the identity $f(1/(3^i-1))=2^i/(3^i+2^{i-1})$, which we proved in part (ii) of this theorem, we get
\begin{flalign*}
f\left(\frac{1}{3^j-3^i}\right)&=f\left(\left(\frac{1}{3^i}\right)\left(\frac{1}{3^{j-i}-1}\right)\right)\\
&=\left(\frac{2}{3}\right)^if\left(\frac{1}{3^{j-i}-1}\right)\\
&=\left(\frac{2}{3}\right)^i\left(\frac{2^{j-i}}{3^{j-i}+2^{j-i-1}}\right)\textrm{\ for\ }j>i>0.\qedhere
\end{flalign*}
\end{proof}
\end{theorem}

\section{Integral Identities}
In this section, we will study the antiderivative $\displaystyle F(x)=\int_0^xf(t)\,dt$.  Our goal is to find inductive formulas describing $F(x)$.  To do this, we first need to prove a result for $F$ analogous to Theorem 1 for $f$.  Then we will use this result in conjunction to prove results similar to Propositions (1)--(3).

Since the integral $\displaystyle\int_0^xf(t)\,dt$ measures the area under a self-similar curve, it exhibits a degree of self-similarity itself.  It turns out that this is the case: We can derive four identities for the integral from our identities for $f$---one of which corresponds to Theorem 1, and three which correspond to Propositions (1)--(3) and serve as iterative formulas for $F$.

\begin{theorem}
For all $x\in[0,1]$, $\displaystyle\int_x^{1-x}f(t)\,dt=1/2-x$.
\begin{proof}
By Theorem 1, we know that for any $t\in[0,1]$, $f(1-t)=1-f(t)$.  So clearly
\begin{flalign*}
\int_a^bf(1-t)\,dt&=\int_a^b1-f(t)\,dt\\
&=[b-a]-\int_a^bf(t)\,dt
\end{flalign*}
for all $a,b\in[0,1]$. So if we let $a=x$ and $b=1-x$, where $x\in[0,1]$, we get
\begin{flalign*}
\int_x^{1-x}\!\!\!\!\!\!f(1-t)\,dt&=[(1-x)-x]-\int_x^{1-x}\!\!\!\!\!\!f(t)\,dt\\
&=[1-2x]-\int_x^{1-x}\!\!\!\!\!\!f(t)\,dt.
\end{flalign*}
By $u$-substitution on $\int_x^{1-x}f(1-t)\,dt$, we have
\begin{flalign*}
\int_x^{1-x}\!\!\!\!\!\!f(1-t)\,dt&=-\int_{1-(x)}^{1-(1-x)}\!\!\!\!\!\!\!\!\!\!\!\!\!\!\!f(t)\,dt\\
&=-\int_{1-x}^{x}\!\!\!f(t)\,dt\\
&=\int_x^{1-x}\!\!\!\!\!\!f(t)\,dt.
\end{flalign*}
Then by substitution,
\begin{equation*}
\int_x^{1-x}\!\!\!\!\!\!f(t)\,dt=[1-2x]-\int_x^{1-x}\!\!\!\!\!\!f(t)\,dt.
\end{equation*}
So
\begin{equation*}
2\int_x^{1-x}\!\!\!\!\!\!f(t)\,dt=1-2x
\end{equation*}
and thus,
\begin{equation*}
\int_x^{1-x}\!\!\!\!\!\!f(t)\,dt=\frac{1}{2}-x, \textrm{\ for all\ }x\in[0,1].\qedhere
\end{equation*}
\end{proof}
\end{theorem}
\noindent
This theorem is illustrated in Fig. 4.  One notable result immediately follows:
\newtheorem*{corollary}{Corollary}
\begin{corollary}
The area under the graph of Bourbaki's function is
\begin{equation*}
A=\int_0^1f(t)\,dt=\frac{1}{2}.
\end{equation*}
\end{corollary}

\begin{figure}[h]
\begin{center}
\includegraphics[width=.5\textwidth]{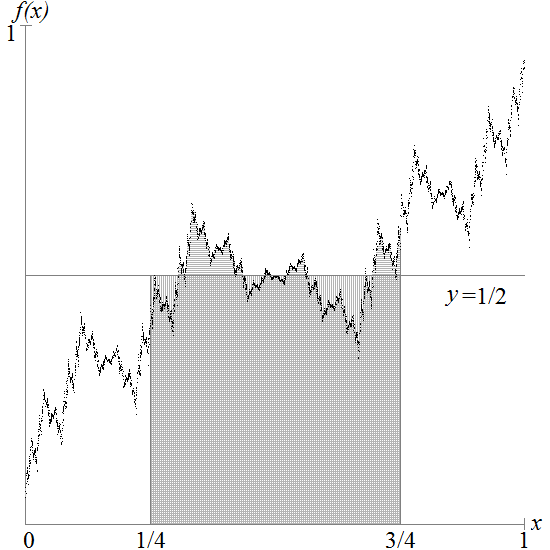}
\caption{The symmetry of the curve generated by $f$ applies to the area under it, as well; over any region $[x,1-x]$ for $x\in[0,1]$, the area under the curve is equal to the area under the line $y=1/2$.}
\end{center}
\end{figure}

\pagebreak

\begin{proposition}
For all $x\in[0,1]$ and $i\geq0$, $\displaystyle\int_0^{x/3^i}\!\!\!\!\!\!f(t)\,dt=\left(\frac{2}{9}\right)^i\int_0^xf(t)\,dt$.
\begin{proof}
By Proposition 1, for all $t,x\in[0,1]$ and $i\geq0$,
\begin{flalign*}
\int_0^xf\left(\frac{t}{3^i}\right)\,dt&=\int_0^x\left(\frac{2}{3}\right)^if(t)\,dt\\
&=\left(\frac{2}{3}\right)^i\int_0^xf(t)\,dt.
\end{flalign*}
Now
\begin{flalign*}
\int_0^xf\left(\frac{t}{3^i}\right)\,dt&=3^i\int_0^{x/3^i}\!\!\!\!\!\!f(t)\,dt\\
&=\left(\frac{2}{3}\right)^i\int_0^xf(t)\,dt
\end{flalign*}
and thus,
\begin{equation*}
\int_0^{x/3^i}\!\!\!\!\!\!f(t)\,dt=\left(\frac{2}{9}\right)^i\int_0^xf(t)\,dt.\qedhere
\end{equation*}
\end{proof}
\end{proposition}
This result is illustrated in Fig. 5.

\begin{figure}[h]
\begin{center}
\includegraphics[width=1.0\textwidth]{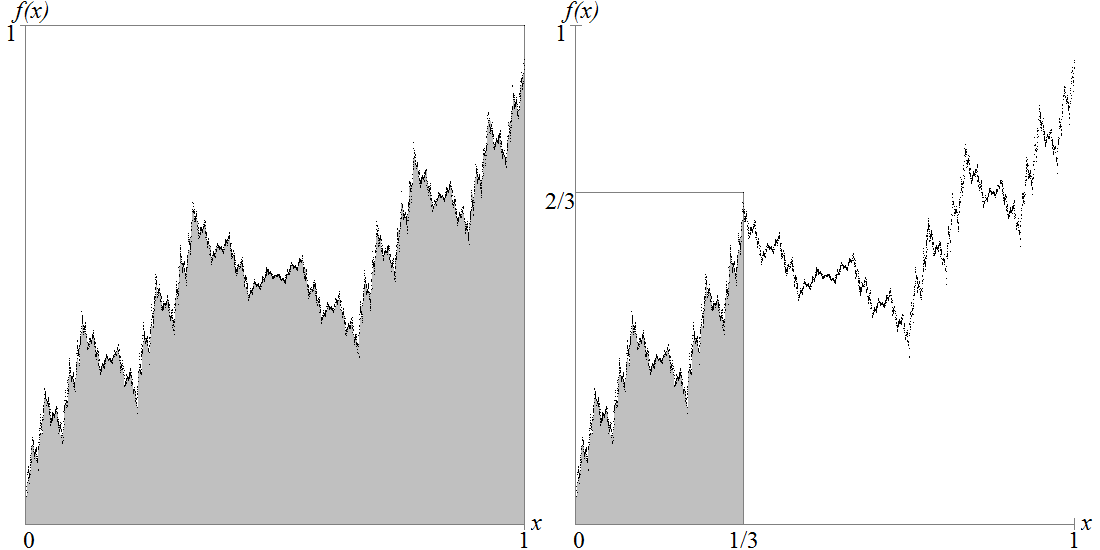}
\caption{Proposition 4 for $i=1$.  The area under the graph over $[0,1]$ and the area over $[0,1/3]$ are in the proportion 1:2/9.  This is exactly the proportion of the area in the two boxes pictured.}
\end{center}
\end{figure}

\pagebreak

\begin{proposition}
For all $x\in[0,1]$ and $i>0$, 
\begin{equation*}
\int_{(2-x)/3^i}^{2/3^i}\!\!\!\!\!\!\!\!\!\!\!\!f(t)\,dt=\frac{2^{i-1}}{9^i}\left[x+\int_0^xf(t)\,dt\right].
\end{equation*}
\begin{proof}
We know by Proposition 2 that for all $t\in[0,1]$ and $i>0$, $f([2-t]/3^i)=(2^{i-1}/3^i)[1+f(t)]$.  So clearly
\begin{flalign*}
\int_0^xf\left(\frac{2-t}{3^i}\right)\,dt&=\int_0^x\frac{2^{i-1}}{3^i}[1+f(t)]\,dt\\
&=\frac{2^{i-1}}{3^i}\int_0^x1+f(t)\,dt\\
&=\frac{2^{i-1}}{3^i}\left[x+\int_0^xf(t)\,dt\right]
\end{flalign*}
for all $x\in[0,1]$. By $u$-substitution on $\int_0^xf([2-t]/3^i)\,dt$, we have
\begin{flalign*}
\int_0^xf\left(\frac{2-t}{3^i}\right)\,dt&=-3^i\int_{2/3^i}^{(2-x)/3^i}\!\!\!\!\!\!\!\!\!\!\!\!\!\!\!f(t)\,dt\\
&=3^i\int_{(2-x)/3^i}^{2/3^i}\!\!\!\!\!\!\!\!\!\!\!\!f(t)\,dt.
\end{flalign*}
So 
\begin{equation*}
3^i\int_{(2-x)/3^i}^{2/3^i}\!\!\!\!\!\!\!\!\!\!\!\!f(t)\,dt=\frac{2^{i-1}}{3^i}\left[x+\int_0^xf(t)\,dt\right]
\end{equation*}
and thus,
\begin{equation*}
\int_{(2-x)/3^i}^{2/3^i}\!\!\!\!\!\!\!\!\!\!\!\!f(t)\,dt=\frac{2^{i-1}}{9^i}\left[x+\int_0^xf(t)\,dt\right],  \textrm{\ for all\ }x\in[0,1] \textrm{\ and\ }i\geq0.\qedhere
\end{equation*}
\end{proof}
\end{proposition}

\begin{proposition}
For all $x\in[0,1]$ and $i>0$, 
\begin{equation*}
\int_{(2+x)/3^i}^{2/3^i}\!\!\!\!\!\!\!\!\!\!\!\!f(t)\,dt=\left(\frac{2^{i-1}}{9^i}\right)x+\left(\frac{2}{9}\right)^i\int_0^xf(t)\,dt.
\end{equation*}
\begin{proof}
This can be proven in the same manner as Proposition 5 if we apply Proposition 3 instead of Proposition 2.
\end{proof}
\end{proposition}
Using Propositions 4--6, we can construct a simple inductive formula for the antiderivative of $f$.

\pagebreak

\begin{theorem}
The antiderivative $F$ of $f$ can be expressed as $F(x)=\displaystyle\lim_{i\to\infty}F_i(x)$, where $F_i$ is defined at any iteration $i\geq0$ as follows: $F_0(x)=x/2$ for all $x\in[0,1]$, every $F_i$ is continuous on $[0,1]$, every $F_i$ is affine on each subinterval $[k/3^i,(k+1)/3^i]$ where $k\in\{0,1,2,\dots,3^i-1\}$, and
\begin{flalign}
F_{i+1}\left(\frac{k}{3^i}\right)&=F_i\left(\frac{k}{3^i}\right),\\
F_{i+1}\left(\frac{k/3^i}{3}\right)&=\frac{2}{9}F_i\left(\frac{k}{3^i}\right),\\
F_{i+1}\left(\frac{1+k/3^i}{3}\right)&=\frac{1}{9}\left[1+\frac{2k}{3^i}-F_i\left(\frac{k}{3^i}\right)\right],\\
F_{i+1}\left(\frac{2+k/3^i}{3}\right)&=\frac{1}{9}\left(\frac{5}{2}+\frac{k}{3^i}\right)+\frac{2}{9}F_i\left(\frac{k}{3^i}\right)\\
F_{i+1}\left(\frac{k+1}{3^i}\right)&=F_i\left(\frac{k+1}{3^i}\right)
\end{flalign}
\begin{proof}
Given the domain of $f$, we will let the antiderivative $F(x)=\displaystyle\int_0^xf(t)\,dt$.  Using this notation, Proposition 4 can be expressed as $\displaystyle F\left(\frac{x}{3^i}\right)=\left(\frac{2}{9}\right)^iF(x)$.  We also can rewrite Propositions 5 and 6 accordingly, but we must adjust them so that their integrals have a lower bound of 0.  
We can do this easily for Proposition 5 using a few substitutions:
\begin{flalign*}
F\left(\frac{1+x}{3^i}\right)&=\int_0^{1/3^i}\!\!\!\!\!\!f(t)\,dt+\int_{1/3^i}^{(1+x)/3^i}\!\!\!\!\!\!\!\!\!\!\!\!f(t)\,dt\\
&=\int_0^{1/3^i}\!\!\!\!\!\!f(t)\,dt+\int_{1/3^i}^{(2-[1-x])/3^i}\!\!\!\!\!\!\!\!\!\!\!\!f(t)\,dt\\
&=\int_0^{1/3^i}\!\!\!\!\!\!f(t)\,dt+\int_{1/3^i}^{2/3^i}\!\!\!\!\!\!f(t)\,dt-\int_{(2-[1-x])/3^i}^{2/3^i}\!\!\!\!\!\!\!\!\!\!\!\!f(t)\,dt\\
&=\left(\frac{2}{9}\right)^i\left(\frac{1}{2}\right)+\frac{2^{i-1}}{9^i}\left(\frac{3}{2}\right)-\frac{2^{i-1}}{9^i}[(1-x)+F(1-x)]\\
&=\frac{2^{i-1}}{9^i}\left[\frac{3}{2}+x-F(x)-\frac{1}{2}+x\right]\\
&=\frac{2^{i-1}}{9^i}[1+2x-F(x)].
\end{flalign*}
Working out Proposition 6 is even simpler:
\begin{flalign*}
F\left(\frac{2+x}{3^i}\right)&=\int_0^{(2+x)/3^i}\!\!\!\!\!\!\!\!\!\!\!\!f(t)\,dt\\
&=\int_0^{1/3^i}\!\!\!\!\!\!f(t)\,dt+\int_{1/3^i}^{2/3^i}\!\!\!\!\!\!f(t)\,dt+\int_{2/3^i}^{(2+x)/3^i}\!\!\!\!\!\!\!\!\!\!\!\!f(t)\,dt\\
&=\left(\frac{2}{9}\right)^i\left(\frac{1}{2}\right)+\frac{2^{i-1}}{9^i}\left(\frac{3}{2}\right)+\frac{2^{i-1}}{9^i}x+\left(\frac{2}{9}\right)^iF(x)\\
&=\frac{2^{i-1}}{9^i}\left(\frac{5}{2}+x\right)+\left(\frac{2}{9}\right)^iF(x).
\end{flalign*}
For $i=1$, the expressions in Propositions 4--6 describe the area under the graph of $f$ over each third of $[0,1]$ in terms of the area over $[0,1]$.  For $i=2$, Propositions 4--6 can be applied over one another to describe the area under the graph over each third of each third of $[0,1]$ in terms of the areas for $i=1$, and so on. Using our three rewritten propositions, we can approximate the graph of $F$ with continuous, affine iterations. We start with the area over $[0,1]$: We know that $F(0)=0$, and from the corollary to Theorem 3, we have $F(1)=1/2$, so our first iteration must be the graph of $F_0(x)=x/2$.  By applying Propositions 4--6 from here, we can evaluate $F_i(k/3^i)$ for any $k\in\{0,1,2,\ldots,3^i-1\}$.  Finally, because the set of all $k/3^i$ is dense in $[0,1]$ as $i$ goes to infinity, we have $\displaystyle F(x)=\lim_{i\to\infty}F_i(x)$ for all $x\in[0,1]$.
\end{proof}
\end{theorem}

Using this theorem, we can obtain a decent approximation of the graph of $F$ (see Fig. 6).

\begin{figure}[h]
\begin{center}
\includegraphics[width=1.0\textwidth]{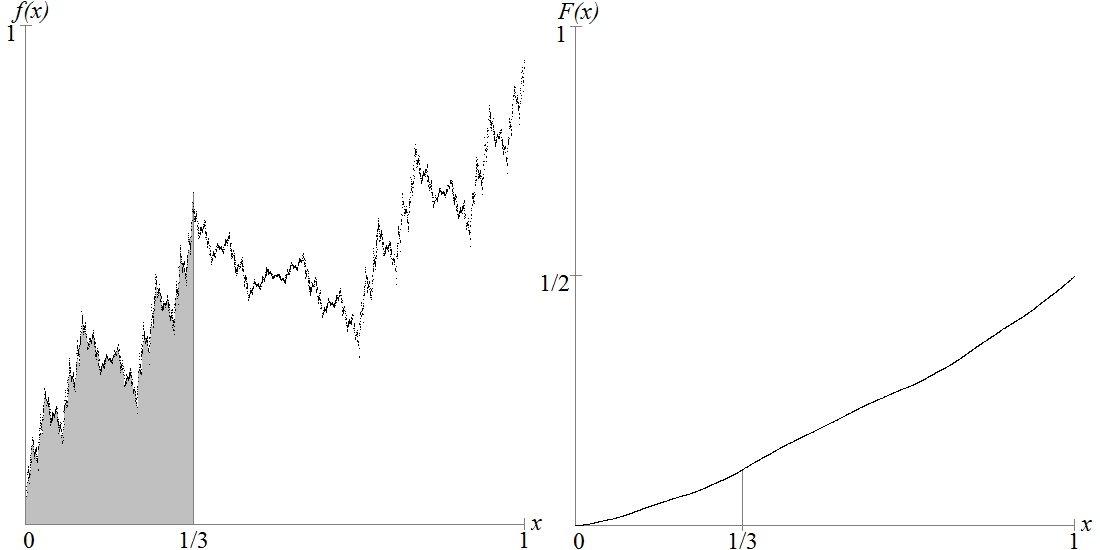}
\caption{$F(x)$ corresponds to the area under the graph of $f$ from 0 to $x$, or $\int_0^xf(t)\,dt$.}
\end{center}
\end{figure}

We can see from the graph that $F$ appears nondecreasing everywhere on $[0,1]$. In fact, this is the case, since $f(x)\geq0$ for all $x\in[0,1]$.  We also can see that the graph looks perfectly smooth, but it also appears to shift between upwards and downwards concavity everywhere.  The Fundamental Theorem of Calculus explains both of these observations.  Okamoto \cite{Okamoto05} has proven that $f$ is continuous and well-defined everywhere on $[0,1]$, and according to the Fundamental Theorem of Calculus, that means that its antiderivative $F$ is continuous, well-defined, and differentiable everywhere on $[0,1]$; also, $F'(x)=f(x)$, so it follows that $F''(x)=f'(x)$.  But Okamoto \cite{Okamoto05} has shown that $f'(x)$ does not exist for any $x\in[0,1]$, so $F''(x)$ also does not exist for any $x\in[0,1]$---in other words, the graph of $F$ is neither concave up nor concave down anywhere on $[0,1]$.  This differs from the concavity of a line, as any linear function of the form $l(x)=ax+b$ will have a second derivative of $l''(x)=0$ and thus could be said to be \emph{both} concave up and concave down.

\pagebreak

\section{Integral Values}
Like $f$, $F$ has no explicit formula, so we must use our identities to predict different values of $F(x)$.

\begin{theorem}
For $i>0$,
\renewcommand{\labelenumi}{(\roman{enumi})}
\begin{enumerate}
\item $\displaystyle\int_0^{1/(3^i+1)}\!\!\!\!\!\!\!\!\!\!\!\!\!\!\!f(t)\,dt=\frac{2^{i-1}}{9^i}\ \frac{3^i-1}{3^i+1}\ \frac{1}{1-(2/9)^i}$,
\item $\displaystyle\int_0^{1/(3^i-1)}\!\!\!\!\!\!\!\!\!\!\!\!\!\!\!f(t)\,dt=\frac{2^{i-1}}{9^i}\ \frac{3^i+1}{3^i-1}\ \frac{1}{1+2^{i-1}/9^i}$,
\item $\displaystyle\int_0^{2/(3^i+1)}\!\!\!\!\!\!\!\!\!\!\!\!\!\!\!f(t)\,dt=\frac{2^{i-1}}{9^i}\ \frac{5\cdot3^i+1}{2\cdot3^i+2}\ \frac{1}{1+2^{i-1}/9^i}$, and
\item $\displaystyle\int_0^{2/(3^i-1)}\!\!\!\!\!\!\!\!\!\!\!\!\!\!\!f(t)\,dt=\frac{2^{i-1}}{9^i}\ \frac{5\cdot3^i-1}{2\cdot3^i-2}\ \frac{1}{1-(2/9)^i}$.
\end{enumerate}
\begin{proof} Let $i>0$.\\

(i) Now, we know by Proposition 4 that for all $x\in[0,1]$ and $i\geq0$,
\begin{equation*}
\int_0^{x/3^i}\!\!\!\!\!\!f(t)\,dt=\left(\frac{2}{9}\right)^i\int_0^xf(t)\,dt
\end{equation*}
and clearly
\begin{equation*}
1-\frac{1}{3^i+1}=\frac{3^i}{3^i+1}
\end{equation*}
So keeping this in mind and applying Theorem 3, we have
\begin{flalign*}
\int_0^{1/(3^i+1)}\!\!\!\!\!\!\!\!\!\!\!\!\!\!\!f(t)\,dt&=\left(\frac{2}{9}\right)^i\int_0^{3^i/(3^i+1)}\!\!\!\!\!\!\!\!\!\!\!\!\!\!\!f(t)\,dt\\
&=\left(\frac{2}{9}\right)^i\left[\int_0^{1/(3^i+1)}\!\!\!\!\!\!\!\!\!\!\!\!\!\!\!f(t)\,dt+\int_{1/(3^i+1)}^{3^i/(3^i+1)}\!\!\!\!\!\!\!\!\!\!\!\!\!\!\!f(t)\,dt\right]\\
&=\left(\frac{2}{9}\right)^i\int_0^{1/(3^i+1)}\!\!\!\!\!\!\!\!\!\!\!\!\!\!\!f(t)\,dt+\left(\frac{2}{9}\right)^i\left(\frac{1}{2}-\frac{1}{3^i+1}\right)
\end{flalign*}
which means that
\begin{equation*}
\left[1-\left(\frac{2}{9}\right)^i\right]\int_0^{1/(3^i+1)}\!\!\!\!\!\!\!\!\!\!\!\!\!\!\!f(t)\,dt=\left(\frac{2^{i-1}}{9^i}\right)\left(\frac{3^i-1}{3^i+1}\right).
\end{equation*}
Since $1-(2/9)^i=0$ when $i=0$, we make the restriction $i>0$ in order to divide on both sides.  This gives us
\begin{equation*}
\int_0^{1/(3^i+1)}\!\!\!\!\!\!\!\!\!\!\!\!\!\!\!f(t)\,dt=\frac{2^{i-1}}{9^i}\ \frac{3^i-1}{3^i+1}\ \frac{1}{1-(2/9)^i},\textrm{\ for\ }i>0.\\
\end{equation*}

(ii) Now, we know by Proposition 5 that for all $x\in[0,1]$ and $i>0$,
\begin{equation*}
\int_{(2-x)/3^i}^{2/3^i}\!\!\!\!\!\!\!\!\!\!\!\!f(t)\,dt=\frac{2^{i-1}}{9^i}\left[x+\int_0^xf(t)\,dt\right]
\end{equation*}
and clearly 
\begin{equation*}
1-\frac{1}{3^i-1}=\frac{3^i-2}{3^i-1}.
\end{equation*}
Thus, we can see that
\begin{equation*}
\int_{1/(3^i-1)}^{2/3^i}\!\!\!\!\!\!\!\!\!\!\!\!f(t)\,dt=\frac{2^{i-1}}{9^i}\left[\frac{3^i-2}{3^i-1}+\int_0^{(3^i-2)/(3^i-1)}\!\!\!\!\!\!\!\!\!\!\!\!f(t)\,dt\right],
\end{equation*}
and so
\begin{equation*}
\int_0^{2/3^i}\!\!\!\!\!\!f(t)\,dt-\int_0^{1/(3^i-1)}\!\!\!\!\!\!\!\!\!\!\!\!f(t)\,dt=\frac{2^{i-1}}{9^i}\left[\frac{3^i-2}{3^i-1}+\int_0^{1/(3^i-1)}\!\!\!\!\!\!\!\!\!\!\!\!f(t)\,dt+\int_{1/(3^i-1)}^{(3^i-2)/(3^i-1)}\!\!\!\!\!\!\!\!\!\!\!\!f(t)\,dt\right].
\end{equation*}
Then
\begin{flalign*}
\int_0^{1/3^i}\!\!\!\!\!\!f(t)\,dt+\int_{1/3^i}^{2/3^i}\!\!\!\!\!\!f(t)\,dt-\int_0^{1/(3^i-1)}\!\!\!\!\!\!\!\!\!\!\!\!f(t)\,dt&=\frac{2^{i-1}}{9^i}\left(\frac{3^i-2}{3^i-1}\right)\\
&+\frac{2^{i-1}}{9^i}\int_0^{1/(3^i-1)}\!\!\!\!\!\!\!\!\!\!\!\!f(t)\,dt\\
&+\frac{2^{i-1}}{9^i}\left(\frac{1}{2}-\frac{1}{3^i-1}\right).
\end{flalign*}
Propositions 4 and 5 give us
\begin{flalign*}
\left(\frac{1}{2}\right)\left(\frac{2}{9}\right)^i+\frac{2^{i-1}}{9^i}\left(1+\frac{1}{2}\right)-\int_0^{1/(3^i-1)}\!\!\!\!\!\!\!\!\!\!\!\!f(t)\,dt&=\frac{2^{i-1}}{9^i}\left(\frac{3^i-2}{3^i-1}\right)\\
&+\frac{2^{i-1}}{9^i}\int_0^{1/(3^i-1)}\!\!\!\!\!\!\!\!\!\!\!\!f(t)\,dt\\
&+\frac{2^{i-1}}{9^i}\left(\frac{1}{2}-\frac{1}{3^i-1}\right).
\end{flalign*}
So
\begin{flalign*}
\left(1+\frac{2^{i-1}}{9^i}\right)\int_0^{1/(3^i-1)}\!\!\!\!\!\!\!\!\!\!\!\!f(t)\,dt&=\left(\frac{1}{2}\right)\left(\frac{2}{9}\right)^i+\frac{2^{i-1}}{9^i}\left(1+\frac{1}{2}\right)-\frac{2^{i-1}}{9^i}\left(\frac{3^i-2}{3^i-1}\right)\\
&-\frac{2^{i-1}}{9^i}\left(\frac{1}{2}-\frac{1}{3^i-1}\right)\\
&=\left(\frac{1}{2}\right)\left(\frac{2}{9}\right)^i+\frac{2^{i-1}}{9^i}\left(1-\frac{3^i-2}{3^i-1}+\frac{1}{3^i-1}\right)\\
&=\left(\frac{1}{2}\right)\left(\frac{2}{9}\right)^i+\left(\frac{2}{9}\right)^i\left(\frac{1}{3^i-1}\right)\\
&=\frac{2^{i-1}}{9^i}\ \frac{3^i+1}{3^i-1}.
\end{flalign*}
Now,
\begin{equation*}
\int_0^{1/(3^i-1)}\!\!\!\!\!\!\!\!\!\!\!\!f(t)\,dt=\frac{2^{i-1}}{9^i}\ \frac{3^i+1}{3^i-1}\ \frac{1}{1+2^{i-1}/9^i},\textrm{\ for\ }i>0.\\
\end{equation*}

(iii) We know by Proposition 5 that for all $x\in[0,1]$ and $i>0$,
\begin{equation*}
\int_{(2-x)/3^i}^{2/3^i}\!\!\!\!\!\!\!\!\!\!\!\!f(t)\,dt=\frac{2^{i-1}}{9^i}\left[x+\int_0^xf(t)\,dt\right]
\end{equation*}
and we can see that clearly
\begin{equation*}
\frac{2-2/(3^i+1)}{3^i}=\frac{2}{3^i+1}.
\end{equation*}
So
\begin{flalign*}
\int_{2/(3^i+1)}^{2/3^i}\!\!\!\!\!\!\!\!\!\!\!\!f(t)\,dt&=\int_0^{2/3^i}\!\!\!\!\!\!f(t)\,dt-\int_0^{2/(3^i+1)}\!\!\!\!\!\!\!\!\!\!\!\!f(t)\,dt\\
&=\int_0^{1/3^i}\!\!\!\!\!\!f(t)\,dt+\int_{1/3^i}^{2/3^i}\!\!\!\!\!\!f(t)\,dt-\int_0^{2/(3^i+1)}\!\!\!\!\!\!\!\!\!\!\!\!f(t)\,dt\\
&=\left(\frac{2}{9}\right)^i\left(\frac{1}{2}\right)+\frac{2^{i-1}}{9^i}\left(1+\frac{1}{2}\right)-\int_0^{2/(3^i+1)}\!\!\!\!\!\!\!\!\!\!\!\!f(t)\,dt
\end{flalign*}
and
\begin{flalign*}
\int_{2/(3^i+1)}^{2/3^i}\!\!\!\!\!\!\!\!\!\!\!\!f(t)\,dt&=\frac{2^{i-1}}{9^i}\left[\frac{2}{3^i+1}+\int_0^{2/(3^i+1)}\!\!\!\!\!\!\!\!\!\!\!\!f(t)\,dt\right]\\
&=\frac{2^{i-1}}{9^i}\left(\frac{2}{3^i+1}\right)+\frac{2^{i-1}}{9^i}\int_0^{2/(3^i+1)}\!\!\!\!\!\!\!\!\!\!\!\!f(t)\,dt.
\end{flalign*}
Thus, by substitution,
\begin{equation*}
\left(\frac{2}{9}\right)^i\left(\frac{1}{2}\right)+\frac{2^{i-1}}{9^i}\left(1+\frac{1}{2}\right)-\int_0^{2/(3^i+1)}\!\!\!\!\!\!\!\!\!\!\!\!f(t)\,dt=\frac{2^{i-1}}{9^i}\left(\frac{2}{3^i+1}\right)+\frac{2^{i-1}}{9^i}\int_0^{2/(3^i+1)}\!\!\!\!\!\!\!\!\!\!\!\!f(t)\,dt,
\end{equation*}
and so
\begin{flalign*}
\left(1+\frac{2^{i-1}}{9^i}\right)\int_0^{2/(3^i+1)}\!\!\!\!\!\!\!\!\!\!\!\!f(t)\,dt&=\left(\frac{2}{9}\right)^i\left(\frac{1}{2}\right)+\frac{2^{i-1}}{9^i}\left(\frac{3}{2}\right)-\frac{2^{i-1}}{9^i}\left(\frac{2}{3^i+1}\right)\\
&=\frac{2^{i-1}}{9^i}\left(1+\frac{3}{2}-\frac{2}{3^i+1}\right)\\
&=\frac{2^{i-1}}{9^i}\ \frac{5\cdot3^i+1}{2\cdot3^i+2}.
\end{flalign*}
Therefore,
\begin{equation*}
\int_0^{2/(3^i+1)}\!\!\!\!\!\!\!\!\!\!\!\!f(t)\,dt=\frac{2^{i-1}}{9^i}\frac{5\cdot3^i+1}{2\cdot3^i+2}\ \frac{1}{1+2^{i-1}/9^i},\textrm{\ for\ }i>0.\\
\end{equation*}

(iv) We know by Proposition 6 that for all $x\in[0,1]$ and $i>0$,
\begin{equation*}
\int_{2/3^i}^{(2+x)/3^i}\!\!\!\!\!\!\!\!\!\!\!\!f(t)\,dt=\frac{2^{i-1}}{9^i}x+\left(\frac{2}{9}\right)^i\int_0^xf(t)\,dt
\end{equation*}
and clearly
\begin{equation*}
\frac{2+2/(3^i-1)}{3^i}=\frac{2}{3^i-1}.
\end{equation*}
The rest of the proof follows steps similar to those in part (iii) of this theorem to give us
\begin{equation*}
\int_0^{2/(3^i-1)}\!\!\!\!\!\!\!\!\!\!\!\!f(t)\,dt=\frac{2^{i-1}}{9^i}\ \frac{5\cdot3^i-1}{2\cdot3^i-2}\ \frac{1}{1-(2/9)^i},\textrm{\ for\ }i>0.\qedhere
\end{equation*}
\end{proof}
\end{theorem}

\section{Comments on Dimension}
We have established that the graphs of $f$ and $F$ both exhibit self-similarity and pathological behavior.  In this section, we will use the self-similarity of the graph of $f$ to prove that it exhibits fractal behavior by having a Hausdorff dimension greater than its topological dimension.   To do this, we will show first that the graph of $f$ has box-counting dimension $\log_3 5$, and then we will use the Mass Distribution Principle (i.e. \cite{Falconer03}) to show that the graph's Hausdorff dimension can be no less than its box-counting dimension.

Although we have shown that the graph of $F$ has no second derivative, we know that $F$ is a $C^1$ function, having a derivative that is continuous everywhere.  Because of this, the graph of $F$ consititutes a rectifiable curve, and so it must have Hausdorff dimension 1.  In this section, we will show that the arc length of the graph lies between $\sqrt{5}/2$ and $3/2$.

\begin{theorem}
If $\Gamma$ is the graph of $f$, then its box-counting dimension $\dim_B(\Gamma)=\log_3 5$.
\begin{proof}
Let $\Gamma$ be the graph of $f$.  We must look at $f$ as $\displaystyle\lim_{i\to\infty}f_i$ here, so we will define $\Gamma_i$ as the graph of $f_i$.  Obviously, $\Gamma=\displaystyle\lim_{i\to\infty}\Gamma_i$.

Now we consider the box-counting dimension $\dim_B(\Gamma_i)$ of $\Gamma_i$.  As the affine pieces of $\Gamma_i$ are defined over intervals of length $1/3^i$, we will count how many boxes of side length $\delta_i=1/3^i$ will cover $\Gamma_i$.  For $\Gamma_0$, the graph of $y=x$ for $x\in[0,1]$, the number $N_{\delta_0}$ of boxes required for the cover is clearly 1; with $\delta_0=1$, we have a box covering the entire graph.  For $\Gamma_1$, we count boxes with $\delta_1=1/3$, and we get $N_{\delta_1}=5$; exactly two boxes cover each of the graph's ``tall'' sides, and one box covers the central portion.  For $\Gamma_2$ with boxes of side length $\delta_2=1/9$, we have $N_{\delta_2}=25$; we count four boxes for each of the four tallest sections, two boxes for each of the four second-tallest sections, and one for the centermost section (see Fig. 7).

\begin{figure}[h]
\begin{center}
\includegraphics[width=1.0\textwidth]{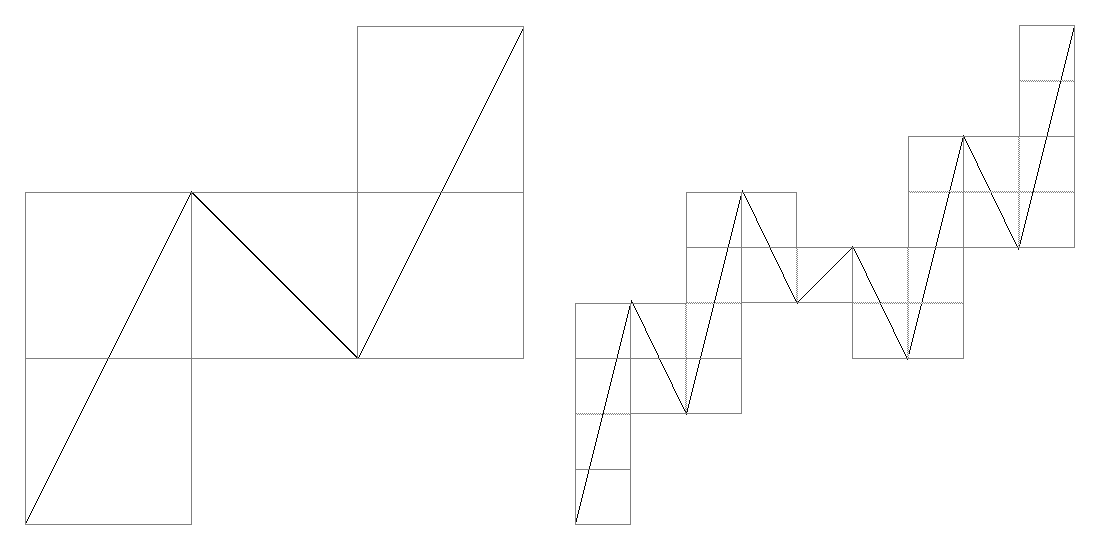}
\caption{Box-counting for $\Gamma_1$ and $\Gamma_2$}
\end{center}
\end{figure}

This gives us an idea of how $i$ varies with $N_{\delta_i}$, but to obtain a general result, we will take the route of Katsuura \cite{Katsuura91} again and view $\Gamma$ as the attractor for a three-component IFS.  We recall that all three contraction mappings from $\Gamma_i$ to $\Gamma_{i+1}$ shinks $\Gamma_i$ horizontally by a factor of $1/3$, but the first and third shrink $\Gamma_i$ vertically by a factor of $2/3$, while the second does so by a factor of only $1/3$.  Now, if we can cover $\Gamma_i$ by $N_{\delta_i}$ boxes, then by necessity, the middle portion of $\Gamma_{i+1}$---the region of the second mapping---could be covered by $N_{\delta_i}$ boxes, as well, since it is $\Gamma_i$ scaled down by a factor of $1/3$ and we are counting how many boxes scaled down by the same factor can cover it.  Applying the similar logic to the regions of the first and third mappings, we can see that $\Gamma_i$ scaled down horizontally by $1/3$ and vertically by $2/3$ will be covered by $2N_{\delta_i}$ boxes scaled down by a factor of $1/3$.  Therefore,
\begin{flalign*}
N_{\delta_{i+1}}&=2N_{\delta_i}+N_{\delta_i}+2N_{\delta_i}\\
&=5N_{\delta_i}
\end{flalign*}
and since $N_{\delta_0}$ is 1, we can say that for $i>0$,
\begin{equation*}
N_{\delta_i}=5^i.
\end{equation*}
Now we consider the formula for box-counting dimension:
\begin{equation*}
\dim_B(\Gamma)=\lim_{\delta_i\to0}\frac{\log(N_{\delta_i})}{-\log(\delta_i)} \textrm{\ (See e.g., \cite{Falconer03}).}
\end{equation*}
And since $\delta=1/3^i$, the formula becomes
\begin{flalign*}
\dim_B(\Gamma)&=\lim_{i\to\infty}\frac{\log(N_{\delta_i})}{-\log(1/3^i)}\\
&=\lim_{i\to\infty}\frac{\log(5^i)}{\log(3^i)}\\
&=\log_3 5.
\end{flalign*}
\end{proof}
\end{theorem}

Now we will prove that the Hausdorff dimension of $\Gamma$ is equal to its box-counting dimension.

\begin{theorem}
If $\Gamma$ is the graph of $f$, then its Hausdorff dimension $\dim_H(\Gamma)=\log_3 5$.
\begin{proof}
To begin, we will consider an alternate iterative construction of $\Gamma$ using Katsuura's mappings.  Let $E_0=[0,1]\times[0,1]$ and define further levels of the construction by $E_{i+1}=w_1(E_i)\cup w_2(E_i)\cup w_3(E_i)$ where $i>0$ and $w_1$, $w_2$, and $w_3$ are the mappings given in \cite{Katsuura91}.  We see that $E_{i+1}\subset E_i$ for all $i\geq0$, and $\displaystyle\bigcap_{i=0}^\infty E_i=\Gamma$ (See Figs. 8 and 9).

\begin{figure}[h]
\begin{center}
\includegraphics[width=1.0\textwidth]{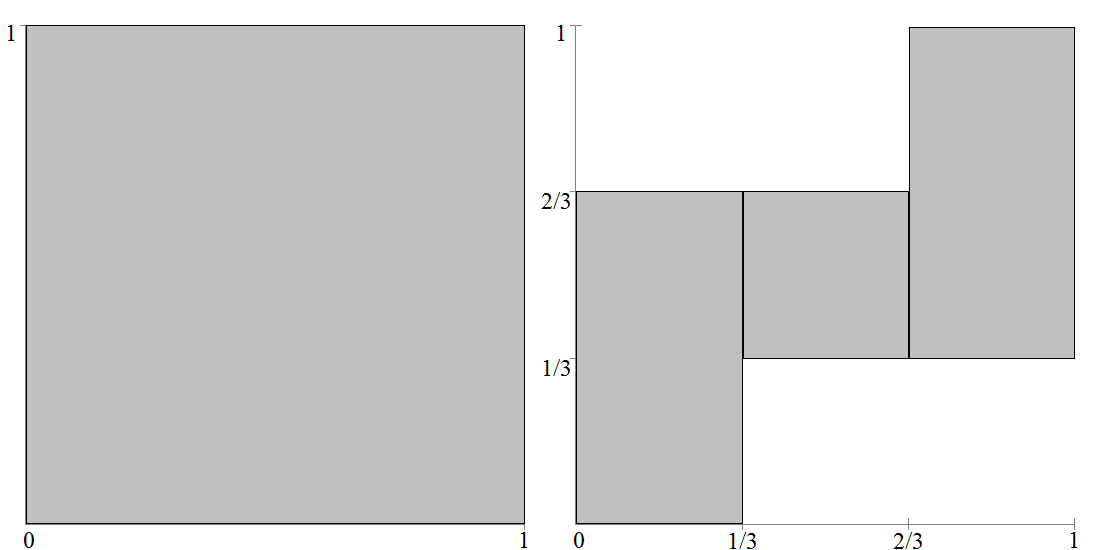}
\caption{$E_0$ and $E_1$.  Note that we can divide the $i$th level of the construction into $3^i$ rectangles of length $(1/3)^i$.}
\end{center}
\end{figure}

\begin{figure}[h]
\begin{center}
\includegraphics[width=1.0\textwidth]{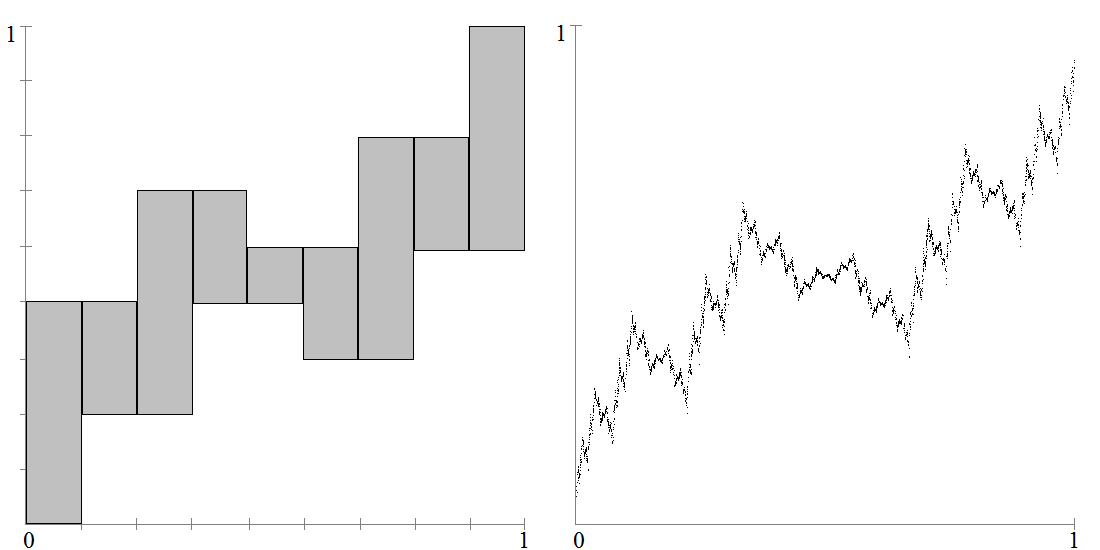}
\caption{$E_2$ and $\Gamma$.  Recall that in Okamoto's construction, linear segments are constructed ``upwards'' to a graph with infinite length, whereas in this construction, rectangular regions are constructed ``downwards'' to a graph with zero area.}
\end{center}
\end{figure}

Using methods related to the box-counting process in Theorem 7, it can be shown that the area of $E_{i+1}$ can be expressed as
\begin{flalign*}
A(E_{i+1})&=\frac{2}{9}A(E_i)+\frac{1}{9}A(E_i)+\frac{2}{9}A(E_i)\\
&=\frac{5}{9}A(E_i),
\end{flalign*}
and since it is obvious that $A(E_0)=1$, we have $A(E_i)=(5/9)^i$ for all $i\geq0$.

Now, let $\mu$ be the natural mass distribution on $\Gamma$; we start with unit mass on $E_0$  and repeatedly ``spead'' this mass over the total area of each $E_i$.  Also, let $U$ be any set whose diameter $|U|<1$.  Then there exists some $i\geq0$ such that
\begin{equation*}
\left(\frac{1}{3}\right)^{i+1}\leq|U|<\left(\frac{1}{3}\right)^i,
\end{equation*}
an inequality that applies to any $U$ satisfying $0<|U|<1$.  From this point, it is clear that for every set $U$ of this type, there is some $i$ such that $U$ is contained in an open square of side length $(1/3)^i$ and $U$ contains points in at most two level-$i$ ``sub-rectangles'' (See Fig. 10).

\begin{figure}[h]
\begin{center}
\includegraphics[width=1.0\textwidth]{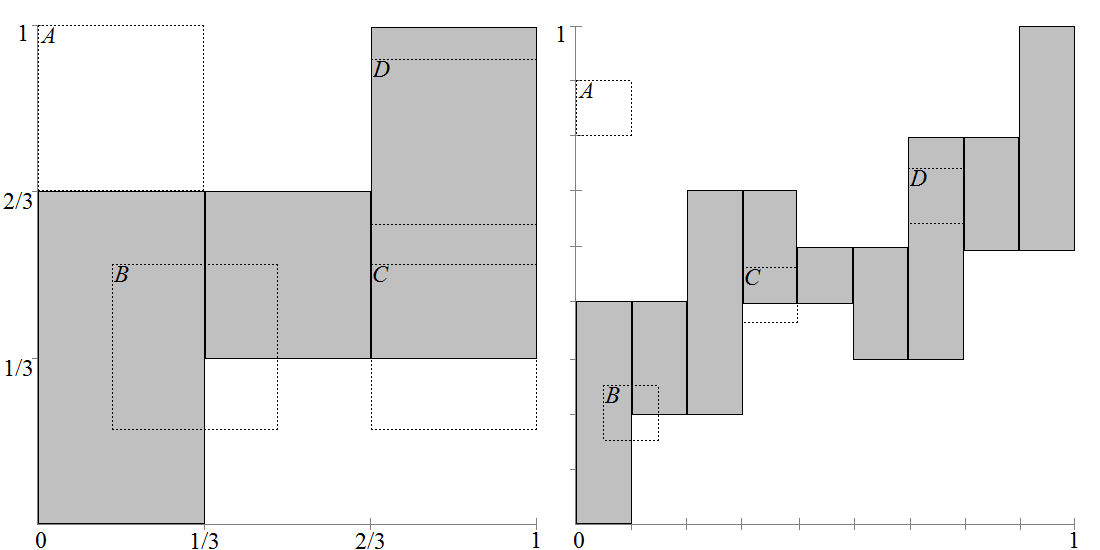}
\caption{Estimating the Hausdorff dimension of $\Gamma$ using the Mass Distribution Principle.  Note that any appropriately-sized set $U$ will ``fit'' in some corresponding open square at some level $E_i$ of the construction, and as a result, $U$ will share points with at most two sub-rectangles of $E_i$.}
\end{center}
\end{figure}

Hence, the area of $U$ is bounded above by the area of the open square containing it; that is, $A(U)\leq(1/9)^i$.  In terms of measure, we know that the entire area of $U$ can be contained in $E_i$, so 
\begin{flalign*}
\mu(U)&\leq\frac{A(U\cap E_i)}{A(E_i)}\\
&\leq\frac{(1/9)^i}{(5/9)^i}\\
&\leq\left(\frac{1}{5}\right)^i.
\end{flalign*}
And since $\displaystyle\left(\frac{1}{3}\right)^{i+1}\leq|U|$ implies that $\displaystyle\left(\frac{1}{3}\right)^i\leq3|U|$, we have
\begin{equation*}
\mu(U)\leq\left(\frac{1}{5}\right)^i=\left(\frac{1}{3^i}\right)^{\log_3 5}\leq(3|U|)^{\log_3 5}=5|U|^{\log_3 5},
\end{equation*}
and thus, by the Mass Distribution Principle, $\log_3 5\leq\dim_H(\Gamma)\leq\dim_B(\Gamma)$, and given the upper bound obtained in Theorem 7, we have $\dim_H(\Gamma)=\log_3 5$.
\end{proof}
\end{theorem}

Because the graph of $F$ is a rectifiable curve, we can determine its arc length over $[0,1]$.  In the following theorem, we will show that this arc length is finite by bounding it above and below.

\begin{theorem}
If $G$ is the graph of $F$ and $L$ is the arc length of $G$, then \linebreak
 $\displaystyle\frac{\sqrt{5}}{2}\leq L\leq\frac{3}{2}$.
\begin{proof}
Let $G$ be the graph of $F$.  We must consider $F$ as the limit of its iterations here, so we will define $G_i$ as the graph of $F_i$.  Obviously, $G=\displaystyle\lim_{i\to\infty}G_i$.

Because each $F_i$ is affine on $[0,1/3^i],[1/3^i,2/3^i],\dots,[(3^i-1)/3^i,1]$ we can apply the Triangle Inequality to the linear ``piece'' of $G_i$ at each of these intervals; for instance, if $l$ is the length of $G_i$ on $[1/3^i,2/3^i]$, we have $|2/3^i-1/3^i|+| F(2/3^i)-F(1/3^i)|\geq l$ (see Fig. 11).

\begin{figure}[h]
\begin{center}
\includegraphics[width=1.0\textwidth]{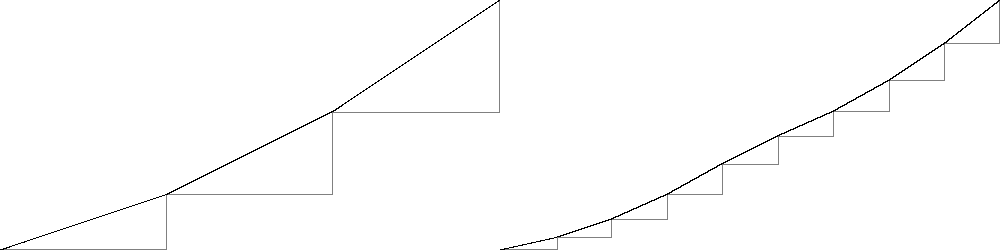}
\caption{Applying the Triangle Inequality to $G_1$ and $G_2$}
\end{center}
\end{figure}

Because the inequality holds over all of $[0,1/3^i],[1/3^i,2/3^i],\dots,[(3^i-1)/3^i,1]$, it also will hold for the sums of the respective sides of each ``triangle''; that is, if $L_i$ is the total arc length of $G_i$, then
\begin{equation*}
L_i\leq\sum_{j=1}^{3^i}\left|\frac{j}{3^i}-\frac{j-1}{3^i}\right|+\sum_{j=1}^{3^i}\left|F_i\left(\frac{j}{3^i}\right)-F_i\left(\frac{j-1}{3^i}\right)\right|.
\end{equation*}
And by taking the limit as $i$ approaches infinity on both sides, we get
\begin{flalign*}
L&\leq\lim_{i\to\infty}\sum_{j=1}^{3^i}\left|\frac{j}{3^i}-\frac{j-1}{3^i}\right|+\sum_{j=1}^{3^i}\left|F_i\left(\frac{j}{3^i}\right)-F_i\left(\frac{j-1}{3^i}\right)\right|\\
&\leq\lim_{i\to\infty}\sum_{j=1}^{3^i}\frac{1}{3^i}+\sum_{j=1}^{3^i}\left|F_i\left(\frac{j}{3^i}\right)-F_i\left(\frac{j-1}{3^i}\right)\right|\\
&\leq\lim_{i\to\infty}1+\sum_{j=1}^{3^i}\left|F_i\left(\frac{j}{3^i}\right)-F_i\left(\frac{j-1}{3^i}\right)\right|.
\end{flalign*}
We observed earlier that $F$ is nondecreasing.  Now, $F_{i+1}(k/3^i)=F_i(k/3^i)$, so by induction, $F(x)=F_i(x)$ wherever $x=k/3^i$.  And since every $F_i$ is affine everywhere between such points, every $F_i$ must be nondecreasing everywhere on $[0,1]$, as well.  This means that for all $j\in\{1,2,\dots,3^i\}$,
\begin{equation*}
F_i\left(\frac{j}{3^i}\right)-F_i\left(\frac{j-1}{3^i}\right)\geq0
\end{equation*}
and thus, we can make the substitution
\begin{equation*}
\left|F_i\left(\frac{j}{3^i}\right)-F_i\left(\frac{j-1}{3^i}\right)\right|=F_i\left(\frac{j}{3^i}\right)-F_i\left(\frac{j-1}{3^i}\right)
\end{equation*}
which gives us
\begin{flalign*}
L&\leq\lim_{i\to\infty}1+\sum_{j=1}^{3^i}F_i\left(\frac{j}{3^i}\right)-F_i\left(\frac{j-1}{3^i}\right)\\
&\leq\lim_{i\to\infty}1+F_i\left(\frac{1}{3^i}\right)-F_i\left(\frac{0}{3^i}\right)+F_i\left(\frac{2}{3^i}\right)-F_i\left(\frac{1}{3^i}\right)\\
&\hspace{3.25 pc}+\dots+F_i\left(\frac{3^i}{3^i}\right)-F_i\left(\frac{3^i-1}{3^i}\right)\\
&\leq\lim_{i\to\infty}1+F_i(1)-F_i(0)\\
&\leq\lim_{i\to\infty}1+\frac{1}{2}-0\\
&\leq\frac{3}{2}.
\end{flalign*}
For the lower bound of $L$, we need only recall that the shortest distance between two points is a straight line.  The endpoints of $G$ are $(0,0)$ and $(1,1/2)$, and the line connecting them has length $\displaystyle\frac{\sqrt{5}}{2}$.  Thus, $L\geq\displaystyle\frac{\sqrt{5}}{2}$.

Now, because the affine segments of each $G_{i+1}$ deviate from the straight lines in $G_i$ from which they are constructed, we see that $L_{i+1}>L_i$ for all $i$.  Approximating $L$, we have $L_2\approx1.1269$, which is strictly greater than $\displaystyle\frac{\sqrt{5}}{2}$.  Thus, $\displaystyle\frac{\sqrt{5}}{2}<L\leq\frac{3}{2}$.
\end{proof}
\end{theorem}
We conclude that although the graph of $F$ exhibits self-similarity and pathological behavior, it is by definition not a fractal.  If we see $F$ as a measure of the area bounded by the graph of $f$ and the $x$-axis, this conclusion makes more sense; a region with a fractal boundary of infinite length still can contain a finite area. (A standard example of this phenomenon is the Koch snowflake; see e.g., \cite{Edgar08}.)

\section{Concluding Remarks}
Okamoto \cite{Okamoto05} shows that Bourbaki's Function is just one member of a parametrized family of functions $F_a$ with analogous constructions.  Using generalizations of expressions (1)--(7), it is possible to prove that given $x\in[0,1]$ and $i>0$, $F_a$ abides by the following rules for all $a\in(0,1)$:
\begin{itemize}
\item $F_a(1-x)=1-F_a(x)$\\
\item $\displaystyle F_a\left(\frac{x}{3^i}\right)=a^iF_a(x)$\\
\item $\displaystyle F_a\left(\frac{2-x}{3^i}\right)=(2a^i-a^{i-1})F_a(x)+(a^{i-1}-a^i)$\\
\item $\displaystyle F_a\left(\frac{2+x}{3^i}\right)=a^iF_a(x)+(a^{i-1}-a^i)$\\
\item $\displaystyle F_a\left(\frac{1}{3^i+1}\right)=\frac{a^i}{1+a^i}$\\
\item $\displaystyle F_a\left(\frac{1}{3^i-1}\right)=\frac{a^{i-1}-a^i}{1+2a^i-a^{i-1}}$\\
\item $\displaystyle F_a\left(\frac{2}{3^i+1}\right)=\frac{1-2a^{i-1}+a^i}{a^{i-1}-2a^i}$\\
\item $\displaystyle F_a\left(\frac{2}{3^i-1}\right)=\frac{a^{i-1}-a^i}{1-a^i}$\\
\item $\displaystyle F_a\left(\frac{1}{3^j-3^i}\right)=\frac{a^{i-1}-a^i}{1-a^i}$
\end{itemize}

Similar generalized identities for the antiderivative of any $F_a$ can be derived using methods similar to those used in this paper. It should be noted, however, not every $F_a$ is nowhere differentiable---a fact that influences several properties of the family of functions, including the dimension of their graphs and the nature of their derivatives almost everywhere in $[0,1]$.  We expect more general proofs to shed light on these subjects.

Obviously, the formulas for finding function and integral values in Theorems 2 and 4 do not guarantee results for any number in $[0,1]$ or even any rational number in that interval.  We do not know of any shortcut for finding $f(1/7)$, for instance, since $1/7$ cannot be expressed in terms of $1/(3^i+1)$, $1/(3^i-1)$, $2/(3^i+1)$, $2/(3^i-1)$, $1/(3^j+3^i)$, or $1/(3^j-3^i)$.  We are unsure if a simple algorithm can be found for evaluating $f(1/m)$ for any natural number $m$; while we attempted to do this by parts for $f(1/3m)$, $f(1/[3m-1])$, and $f(1/[3m-2])$, we were unsuccessful.

\section{Acknowledgments}
The approximate graphs of $f$ and $F$ were produced using \emph{Dynamical Grapher for Quadratic Maps} \cite{Basselet}.

Many thanks to Professor Daniel Jackson for providing resources, feedback, and encouragment throughout the duration of this project.

\pagebreak

\end{document}